\documentclass[reqno,11pt]{amsart}
\usepackage{epsf, amsmath, amsfonts, amssymb, amsthm, amscd}

\def\jref#1#2#3#4#5{#1 {\bf #2} (19#3), #4--#5.}

\def\mybibitem#1#2#3#4#5#6#7#8#9{\bibitem[#1]{#2}#3, {\sl #4},
\jref{#5}{#6}{#7}{#8}{#9}}
\def\mybibbook#1#2#3#4#5{\bibitem[#1]{#2}#3, {\sl #4}, #5}
\def\mybibprep#1#2#3#4#5{\bibitem[#1]{#2}#3, {\sl #4}, {\tt Preprint} #5.}

\def\cmp{Commun. Math. Phys.}
\def\jsp{J. Stat. Phys. }

\def\ptrf{Probab. Theory Relat. Fields}

\def\bal{\begin{align}}
\def\be{\begin{equation}}
\newcommand{\req}[1]{\eqref{#1}}

\makeatletter
\def\captionfont@{\footnotesize}
\def\captionheadfont@{\scshape}

\long\def\@makecaption#1#2{%
  \vspace{2mm}
  \setbox\@tempboxa\vbox{\color@setgroup
    \advance\hsize-6pc\noindent
    \captionfont@\captionheadfont@#1\@xp\@ifnotempty\@xp
        {\@cdr#2\@nil}{.\captionfont@\upshape\enspace#2}%
    \unskip\kern-6pc\par
    \global\setbox\@ne\lastbox\color@endgroup}%
  \ifhbox\@ne 
    \setbox\@ne\hbox{\unhbox\@ne\unskip\unskip\unpenalty\unkern}%
  \fi
  \ifdim\wd\@tempboxa=\z@ 
    \setbox\@ne\hbox to\columnwidth{\hss\kern-6pc\box\@ne\hss}%
  \else 
    \setbox\@ne\vbox{\unvbox\@tempboxa\parskip\z@skip
        \noindent\unhbox\@ne\advance\hsize-6pc\par}%
\fi
  \ifnum\@tempcnta<64 
    \addvspace\abovecaptionskip
    \moveright 3pc\box\@ne
  \else 
    \moveright 3pc\box\@ne
    \nobreak
    \vskip\belowcaptionskip
  \fi
\relax
}
\makeatother
\def\writefig#1 #2 #3 {\rlap{\kern #1 truecm
\raise #2 truecm \hbox{#3}}}
\def\figtext#1{\smash{\hbox{#1}}
\vspace{-5mm}}

\def\math#1{\ifmmode
\mathchoice{\mbox{$\displaystyle\rm#1$}}
{\mbox{$\textstyle\rm#1$}}
{\mbox{$\scriptstyle\rm#1$}}
{\mbox{$\scriptscriptstyle\rm#1$}}\else
{\mbox{$\rm#1$}}\fi}		

\newcommand{\nnm}{\nonumber}

\DeclareMathSymbol{\leqslant}{\mathalpha}{AMSa}{"36}
\DeclareMathSymbol{\geqslant}{\mathalpha}{AMSa}{"3E}
\DeclareMathSymbol{\doteqdot}{\mathalpha}{AMSa}{"2B}
\DeclareMathSymbol{\circlearrowright}{\mathalpha}{AMSa}{"08}
\DeclareMathSymbol{\eset}{\mathalpha}{AMSb}{"3F}

\renewcommand{\emptyset}{\eset}
\renewcommand{\subset}{\subseteq}

\renewcommand{\leq}{\;\leqslant\;}
\renewcommand{\geq}{\;\geqslant\;}
\newcommand{\defby}{\;\doteqdot\;}

\newcommand{\dd}{\,\text{\rm d}}
\newcommand{\e}[1]{\,\text{\rm e}^{#1}\,}

\newcommand{\inftwo}[2]{\inf_{\substack{#1 \\ #2}}}
\newcommand{\sumtwo}[2]{\sum_{\substack{#1 \\ #2}}}
\newcommand{\sumthree}[3]{\sum_{\substack{#1 \\ #2 \\ #3}}}

\newcommand{\prodtwo}[2]{\prod_{\substack{#1 \\ #2}}}
\newcommand{\prodthree}[3]{\prod_{\substack{#1 \\ #2 \\ #3}}}

\newcommand{\ra}{\rightarrow}

\newcommand{\bk}[1]{{\ifmmode{\langle#1\rangle}\else${\langle#1\rangle}$\fi}} 
\def\sp#1#2{{\ifmmode{\langle{#1}|{#2}\rangle}		
\else${\langle{#1}|{#2}\rangle}$\fi}}

\newtheorem{lem}{Lemma}[section]
\newtheorem{pro}{Proposition}[section]
\newtheorem{thm}{Theorem}[section]
\newtheorem{cor}{Corollary}[section]

\makeatletter
\def\usedefcounter#1{\@nmbrlisttrue\def\@listctr{#1}}

\def\itemiz{\ifnum \@itemdepth >3 \@toodeep\else \advance\@itemdepth \@ne
        \edef\@itemitem{labelitem\romannumeral\the\@itemdepth}%
        \list{\csname\@itemitem\endcsname}{
        \setlength{\topsep}{-0mm}
        \setlength{\parsep}{-0.5mm}
        \setlength{\itemsep}{0mm}
        \setlength{\labelsep}{2mm}
        \settowidth{\leftmargin}{M.}
        \addtolength{\leftmargin}{\labelsep}
        \def\makelabel##1{\hss\llap{##1}}}\fi}

\makeatother

\newenvironment{deflist}
{\begin{list}{{\rm (D\arabic{defcounter})}}
{\settowidth{\leftmargin}{(D222)}
\setlength{\labelsep}{2mm}
\addtolength{\leftmargin}{\labelsep}
\settowidth{\labelwidth}{(D222)}
\usedefcounter{defcounter}}}
{\end{list}}

\newcounter{defcounter}

\def\cA{\ensuremath{\mathcal A}}

\def\cC{\ensuremath{\mathcal C}}
\def\cD{\ensuremath{\mathcal D}}

\def\cO{\ensuremath{\mathcal O}}

\def\frK{\ensuremath{\frak K}}

\def\agr{\ensuremath{\alpha} }

\def\dgr{\ensuremath{\delta} }

\def\egr{\ensuremath{\varepsilon} }

\def\lgr{\ensuremath{\lambda} }

\def\Lgr{\ensuremath{\Lambda} }

\newcommand{\bbE}{{\ensuremath{\mathbb E}} }

\newcommand{\bbJ}{{\ensuremath{\mathbb J}} }

\newcommand{\bbP}{{\ensuremath{\mathbb P}} }

\newcommand{\bbR}{{\ensuremath{\mathbb R}} }

\newcommand{\bbZ}{{\ensuremath{\mathbb Z}} }
\newcommand{\bbZZ}{{\ensuremath{\mathbb Z^2}}}

\newcommand{\normI}[1]{\lVert#1\rVert_{\scriptscriptstyle 1}}
\newcommand{\norminf}[1]{\lVert#1\rVert_{\scriptscriptstyle \infty}}
\newcommand{\abs}[1]{\lvert#1\rvert}

\newcommand{\setof}[2]{\{#1\,:\,#2\}}

\newcommand{\ol}[1]{\overline{#1}}
\newcommand{\ul}[1]{\underline{#1}}

\def\ds{\displaystyle}

\newcommand{\var}{{\rm var}}

\newcommand{\Prob}{{\rm Prob}}

\def\thebox{{\Lgr_L}}

\def\Cpl{\Psi}
\def\pot{V}
\newcommand\cf[1]{\chi(#1)}
\newcommand\rnb[1]{\bk{#1}_{\! r}}

\def\given{\,\big\vert\,}
\def\nea{{\!\scriptscriptstyle \nearrow}}

\def\back{\hspace*{-2mm}}

\newcommand{\ERW}[2]{{\bbE}_{#1}[#2]} 
 
\newcommand{\ERWNG}[3]{{\bbE}_{#1}^{{#2}}[#3]} 
\newcommand{\Acomp}{A^{\rm c}} 
 
\newcommand{\diam}{{\rm diam}}
\newcommand{\sprod}[2]{(#1,#2)}
\newcommand{\bext}{\partial^{\rm ext}}
\newcommand{\AcConn}{{\ds \mathop{\longleftrightarrow}^{\phantom{{}^{\rm c}}
A^{\rm c}}}}
\newcommand{\BConn}{{\ds \mathop{\longleftrightarrow}^B}}
\begin{document}
\title[\,]{Non-Gaussian Surface Pinned by a Weak Potential}
\author{J.-D. Deuschel}
\address{
Fachbereich Mathematik, Sekr. MA 7-4, TU-Berlin, Stra\ss e des 17. Juni 136,
D-10623 Berlin, Germany
}
\email{deuschel\@@math.tu-berlin.de}
\author{Yvan Velenik}
\address{
Fachbereich Mathematik, Sekr. MA 7-4, TU-Berlin, Stra\ss e des 17. Juni 136,
D-10623 Berlin, Germany
}
\email{velenik\@@math.tu-berlin.de}
\date{\today}
\begin{abstract}
We consider a model of a two-dimensional interface of the SOS type, with
finite-range, even, strictly convex, twice continuously differentiable
interactions. We prove that, under an arbitrarily weak potential favouring
zero-height, the surface has finite mean square heights. We consider the cases
of both square well and $\delta$ potentials. These results extend previous
results for the case of nearest-neighbours Gaussian interactions in \cite{DMRR}
and \cite{BB}. We also obtain estimates on the tail of the height distribution
implying, for example, existence of exponential moments. In the case of the
$\dgr$ potential, we prove a spectral gap estimate for linear functionals. We
finally prove exponential decay of the two-point function (1) for strong
$\dgr$-pinning and the above interactions, and (2) for arbitrarily weak
$\dgr$-pinning, but with finite-range Gaussian interactions.
\end{abstract}
\maketitle

\section{Introduction}
Even though the understanding of phase separation and related interfacial
phenomena for two-dimensional systems such as the Ising model, has greatly
improved  recently, the situation for three-dimensional systems remains quite
unsatisfactory. For example, even in the three-dimensional Ising model, several
basic questions remain open: Existence of a roughening transition, proof that
the wetting transition occurs at a non-trivial value of the boundary magnetic
field (or proof of the contrary), or even instability of the
$(1,1,1)$-interface. To gain some insights in these problems, it is very useful
to consider simpler SOS-type, effective models for interfaces. In these models,
the interface is modelized as a function $\ul h$ from a subset of $\bbZ^d$ to
$\bbR$ (or $\bbZ$, but we restrict our attention to the former case) where
$h_i\equiv\ul{h}(i)$ represents the height of the surface above, or below, the
site $i$; the energy associated to this surface is specified by some function
of these heights, $H(\ul h) = \sum_{i,j}\Cpl_{ij}(h_i-h_j)$. Unfortunately,
even these much simplified models remain rather difficult to handle, and most
of the results which have been obtained are restricted to the Harmonic case,
where $\Cpl_{ij}(h_i-h_j)=\tfrac12(h_i-h_j)^2$. It is therefore valuable to
find ways to extend such results to a larger class of models, by providing
arguments which are less sensitive to the particular features of the underlying
interaction.

The aim of the present work is to extend to a large class of interactions
results about the pinning of an interface by a weak potential. We emphasize
that these results do not follow from perturbation around the Gaussian case.
Let $\thebox=[-L,L]^2\cap\bbZZ$. It is well known that the mean square height
at any fixed site $i\in\thebox$ w.r.t. the Gaussian measure with $0$-boundary
condition has a logarithmic divergence in the thermodynamic limit. However, it
was shown in \cite{DMRR} that the addition of an arbitrarily weak
self-potential favouring height $0$ would localize the surface, in the sense
that this mean square height remains bounded uniformly in $L$. We show that
this holds for a large class of finite-range non-Gaussian interactions.

As a byproduct of our technique, we also obtain an improved version of a result
of \cite{BB} proving the existence of a massgap for the model with
nearest-neighbours, Gaussian interactions, in the presence of an arbitrarily
weak pinning potential. Even though we are not able to extend their result to
non-Gaussian interactions, we show that it holds in the following situations:
(1) For a large class of finite-range non-Gaussian interactions but with a
sufficiently strong $\dgr$-pinning, and (2) for a large class of finite-range
Gaussian interactions with arbitrarily weak $\dgr$-pinning.

We finally prove the following new result (as far as we know) on the tail of
the distribution, valid for the same class of finite-range non-Gaussian
interactions submitted either to square-well or $\dgr$ potentials: the
probability that the height of the interface above some site $i$ is larger than
$T$ (large) is bounded from above by $\exp(-\cO(T^2/\log T))$; this implies of
course existence of all moments, including exponential ones. If the interaction
has bounded second derivatives, then we also prove the corresponding lower
bound.

\bigskip
We restrict our attention to dimension two since it is the relevant case to
describe an interface in a three-dimensional medium. It is also the most
interesting one as far as pinning is concerned. Indeed, in dimensions greater
than two, the situation is completely different: The mean square of the height
of the interface is already finite {\em without} a pinning potential. The
behaviour of the two-point function is also different: without pinning potential,
it has a power-law decay. However, the addition of such a potential
would make this decay exponential.

\bigskip
In Section \ref{s_statements}, we define the models and state the main results
of this paper. Proofs of these statements are given in Sections \ref{s_pinV} and
\ref{s_expdecay}. Our main estimate, Proposition \ref{pro_drydecay}, is proved in
Section \ref{s_mainestimate}. Some technical estimates are given in the appendix.
\section{Models and results}\label{s_statements}
Let $r$ be some strictly positive integer, the {\em range} of the interaction.
The interaction between sites $i$ and $j$, $\Cpl_{i,j}(h_j-h_i)$ is supposed to
satisfy the following conditions:
\begin{itemize}
\item{\bf Translation invariance:} $\Cpl_{i,j}=\Cpl_{0,j-i}\equiv\Cpl_{j-i}$\,.
\item{\bf Finite range:} $\Cpl_k \equiv 0$ if $\normI k > r$.
\item{\bf Symmetry:} $\Cpl_k = \Cpl_{-k}$ and $\Cpl_k(x)=\Cpl_k(-x)$.
\item{\bf Smoothness:} $\Cpl_k$ is twice continuously differentiable.
\item{\bf Irreducibility:} $\Cpl_k$ is convex, i.e. $\Cpl_k''(x)\geq 0$.
Moreover, there exists $c>0$ such that the random walk on $\bbZZ$ with
transition rates $P_c(0,j)$ given by  $1$ if $\Cpl_k''(h)\geq c\;\forall h$,
and $0$ otherwise, is irreducible.
\end{itemize}
All these conditions are natural, and standard in this kind of problem, except
for the last one, which happens to be necessary to be able to use standard result
about the random walk appearing in the random walk representation
(described below) of a related Gaussian model obtained using Brascamp-Lieb
inequality.

\noindent
{\bf Remark:} The hypothese on translation invariance could be removed easily.
We only left it for notational convenience.

\medskip
Let $b\in\bbR$ and let
$\Lgr\Subset\bbZZ$. The Gibbs measure with $b$-b.c.
in $\Lgr$ is the probability measure on $\bbR^\bbZZ$ given by
\begin{equation}
\mu^b_{\Lgr}(\dd \ul h) \defby \frac 1 {Z^b_\Lgr}
\exp\bigl\{-\sum_{\rnb{ij}\cap\Lgr\neq\emptyset} \Cpl_{j-i}(h_i-h_j) \bigr\}
\prod_{i\in\Lgr} \dd h_i \prod_{i\not\in\Lgr}\dgr_{b}(\dd h_i)\,,
\end{equation}
where $\rnb{ij}$ is any pair of distinct sites $i$ and $j$ such that
$\norminf{j-i}\leq r$ and $\dgr_b$ is the point-mass at $b$. Expectation value
and variance with respect to $\mu^b_\Lgr$ are denoted by
$\bk{\,\cdot\,}^b_\Lgr$ and $\var^b_\Lgr(\,\cdot\,)$.

Let $\egr$ and $a$ be two strictly positive real numbers; the potential
$\pot:\bbR\ra\bbR$ is defined by
\begin{equation}
\pot(h) \defby -\egr\cf{\abs{h}\leq a}\,,
\end{equation}
where $\cf{\cA}$ is the indicator function of the event $\cA$. The Gibbs
measure with $0$-b.c. on $\Lgr$ and potential $\pot$ is the
probability measure defined by
\begin{equation}
\mu^\pot_\Lgr(\dd \ul h) \defby \frac 1{Z^\pot_\Lgr}
\exp\{-\sum_{i\in\Lgr}\pot(h_i)\}\;\mu^0_\Lgr(\dd \ul h)\,.
\end{equation}
Expectation value with respect to this measure is written
$\bk{\,\cdot\,}^\pot_\Lgr$.

\medskip
Let $\thebox\defby [-L,L]^2\cap\bbZZ$. Our first result shows that, for any
$\egr$ and $a$, the mean square height of the field is finite, uniformly in $L$.
This generalizes the corresponding result of \cite{DMRR} valid for
nearest-neighbors Gaussian interactions.
\begin{thm}\label{thm_pinV}
There exists a constant $C_1=C_1(a(e^\egr-1),c,r)<\infty$ such that, for any
$i\in\thebox$ and $\forall L$,
$$
\bk{h_i^2}^\pot_\thebox \leq C_1\,.
$$
Moreover, if $a(e^\egr-1)\sqrt c$ is small, then there exists $C_2=C_2(r)>0$ such that
$C_1 \leq 4 a^2+\frac {C_2}c \abs{\log(a(e^\egr-1)\sqrt c)}$.
\end{thm}
Note that the fact that $C_1$ depends on $a$ only through the product $a\sqrt c$
is natural, since otherwise we could improve the result by rescaling the field
$\ul h$.

In fact, using the same techniques, it is possible to obtain a much stronger
statement about pinning of the field, namely existence of exponential moments.
Indeed, this is a consequence of the following estimates on the tail of the
height distribution.
\begin{thm}\label{thm_tail}
There exist $C_3=C_3(a(e^\egr-1),c,r)$ and $T_0=T_0(a(e^\egr-1),c,r) a$ such that,
for all $T>T_0$ with $T\gg a$, and all $L$,
$$
\mu^V_\thebox(h_i\geq T) \leq e^{-C_3\; T^2/ \log T}\,.
$$
Moreover, if $\Cpl_k''(h)\leq \frac 1c$, for all $k$ and $h$, then there exists
a constant $C_4=C_4(a(e^\egr-1),c,r)<\infty$ such that, for all $T>1$ and all
$L$,
$$
\mu^V_\thebox(h_i\geq t) \geq e^{-C_4\; T^2/ \log T}\,.
$$
\end{thm}

\medskip
In \cite{BB}, a statement analogous to Theorem \ref{thm_pinV} was proved,
together with the exponential decay of the 2-point function, for a slightly
different measure to which we will refer as the $\dgr$-pinning. Their measure
corresponds to
\begin{equation}\label{eq_measdelta}
\mu^J_\thebox(\dd \underline h) = \frac1{Z^J_\thebox}
\exp[-\sum_{\rnb{ij}\subset\thebox}\Cpl(h_i-h_j) -
\sumtwo{\rnb{ij}}{i\in\thebox,\,j\not\in\thebox}\Psi(h_i)] \prod_{i\in\thebox}
(\dd h_i + e^J \delta_0(\dd h_i))\,,
\end{equation}
where $J$ is some real parameter. (In fact, they considered the Gaussian case,
with periodic boundary conditions and nearest-neighbors interaction).
Expectation value and variance with respect to $\mu^J_\thebox$ are written
$\bk{\,\cdot\,}^J_\thebox$ and $\var^J_\thebox(\,\cdot\,)$.

The measure in \req{eq_measdelta} can be seen as the limit of the measure
$\mu^V_\thebox$, when $\egr\ra\infty$ with $2(e^\egr-1)a=e^J$ (using, for
example, Lebesgue's Theorem). Since the bounds in Theorems \ref{thm_pinV} and
\ref{thm_tail} only depend on the product $(e^\egr-1)a$, they readily give the
following\footnote{Notice that the result on pinning in \cite{BB} cannot be
deduced from the corresponding statement in \cite{DMRR} since the bound given
in this last work diverges in the limit $\egr\ra\infty$, $2(e^\egr-1)a=e^J$.}
\begin{cor}\label{cor_pinJ}
There exists a constant $C'_1=C'_1(J,c,r)<\infty$ such that, for any $i\in\thebox$
and for all $L$,
$$
\bk{h_i^2}^J_\thebox \leq C'_1\,.
$$
Moreover, if $e^J\sqrt c$ is small, then there exists $C'_2=C'_2(r)>0$ such that
$C'_1 \leq \frac {C'_2}c \abs{\log(e^J\sqrt c)}$.
\end{cor}
\begin{cor}\label{cor_tail}
There exist $C'_3=C'_3(J,c,r)$ and $T_0=T_0(J,c,r)$ such that, for all $T>T_0$ and all
$L$,
$$
\mu^J_\thebox(h_i\geq T) \leq e^{-C'_3\; T^2/ \log T}\,.
$$
Moreover, if $\Cpl_k''(h)\leq \frac 1c$, for all $k$ and $h$, then there exists
a constant $C'_4=C'_4(J,c,r)<\infty$ such that, for all $T>1$ and all $L$,
$$
\mu^J_\thebox(h_i\geq t) \geq e^{-C'_4\; T^2/ \log T}\,.
$$
\end{cor}
It is in fact possible to obtain bounds on more general quantities. Indeed, we
have the following
\begin{thm}\label{thm_sgap}
There exists a constant $C_5=C_5(J,c,r)<\infty$ such that, for any
$\agr:\thebox\ra\bbR$ and all $L$,
$$
\var^J_\thebox(\sprod\agr h) \leq C_5 \sprod{\agr}{\agr}\,,
$$
where $\sprod \agr h = \sum_{i\in\Lgr} \agr(i) h_i$.
\end{thm}
\noindent
{\bf Remark:} This result implies in particular that, for any
$\Lgr'\subset\thebox$,
\begin{equation}
\var^J_\thebox(\frac 1{\abs{\Lgr'}}\sum_{i\in\Lgr'}h_i) = \cO(\abs{\Lgr'}^{-1})\,.
\end{equation}

It does not seem possible with our techniques to prove the same kind of result
for arbitrary functions. However, it is still possible to get the following
\begin{cor}\label{cor_odd}
Let $C^{1,\rm o}_{\rm b}$ be the set of functions $F\in C^1(\bbR^\thebox)$ such
that $\norminf{f_i}\equiv\norminf{\frac\partial{\partial h_i} F} < \infty$,
$i\in\thebox$, and $F$ is odd in each coordinate, i.e. $F(T_i h) = -F(h)$ where
$(T_i h)_j = (1-2\dgr_{i,j})h_j$. Then there exists a constant
$C'_5=C'_5(J,c,r)<\infty$ such that
$$
\var^J_\thebox(F) \leq C'_5 \sum_{i\in\thebox} \norminf{f_i}^2\,.
$$
\end{cor}

\medskip
The results we are able to obtain about exponential decay of the two-point
function are less satisfactory. Since it is not clear how the technique used in
this paper (and taken from \cite{BB}) should be used to prove exponential decay
of the 2-point function in the non-Gaussian case (the corresponding random walk
representation being much more complicated, see below), we have to restrict our
attention to Gaussian interactions. The result we obtain in this case reads
\begin{thm}\label{thm_expdecayGauss}
Suppose that, in addition to the above hypotheses, the interaction is Gaussian,
i.e. $\Cpl_k(x)=c_k x^2$. Then there exists $C_6=C_6(J,c,r)>0$ such that, for
all $L$,
$$
\bk{h_ih_j}^J_\thebox \leq \tfrac 1{C_6} e^{-C_6\normI{j-i}}\,.
$$
Moreover, if $e^J\sqrt c$ is small, then there exists $C_7=C_7(r)>0$ such that
$C_6 \geq C_7 \abs{\log(e^J\sqrt c)}$.
\end{thm}
This improves the result of \cite{BB} since it holds for $0$-b.c. and
finite-range interactions. In the non-Gaussian case, we are only able to prove
exponential decay in the strong pinning regime, i.e. when the parameter $J$ is
sufficiently large. In this case, the following can be proved,
\begin{thm}\label{thm_expdecaylargeJ}
There exist $J_0$ and $C_8=C_8(J,c,r)>0$ such that, for all $J\geq
J_0$, and all $L$,
$$
\bk{h_ih_j}^J_\thebox \leq \tfrac 1{C_8} e^{-C_8\normI{j-i}}\,.
$$
\end{thm}

\bigskip
The basic strategy to obtain these exponential decay results is taken from
\cite{BB}. Our main contribution is Proposition \ref{pro_drydecay} which
replaces an estimate in \cite{BB} the proof of which, based on reflection
positivity and an entropy estimate, limited their analysis to nearest-neighbors
Gaussian interactions. Our method, more robust, is inspired in part by their
entropy estimate.

\medskip
A basic tool for our analysis is the following {\em random walk representation}
of two-point functions (see \cite{DGI}). For any $\Lgr\Subset\bbZZ$ and any
$i,j\in\Lgr$ (it is possible that $i=j$), the following holds
\begin{align}
\bk{h_ih_j}^0_\Lgr = \ERWNG{i}{\Lgr,0}{\int_0^{\tau_\Lgr} \cf{\eta_s=j}\dd s}\,,
\end{align}
where $\eta$ is a random walk on $\bbZZ$ starting at $i$, $\eta_s$ its position
at time $s$ and $\tau_\Lgr = \inf\setof{s\geq 0}{\eta_s\not\in\Lgr}$.
Expectation of an event $\cA$ depending only on $\eta$ is given by
\begin{equation}
\ERWNG{i}{\Lgr,0}{\cA} = \bk{\ERW{i,\cdot}{\cA}}^0_\Lgr\,,
\end{equation}
with $\ERW{i,\ul h}{\,\cdot\,}$ denoting joint expectation w.r.t. the symmetric
diffusion $\ul h(s)$ and the random walk in $\bbZZ$, starting at $i$, with
jump-rate at time $s$ given by
\begin{equation}
p_{\ul h(s)}(i,j) = 2\Cpl_{j-i}''(h_j(s)-h_i(s))\,.
\end{equation}
Observe that in the Gaussian case $\Cpl_k''=c_k$ is independent of $\ul h$ and
therefore $\ERWNG{i}{\Lgr,0}{\cA}\equiv\ERW{i}{\cA}$, the expectation w.r.t. the
random walk in the plane starting at $i$, with jump-rate $p(i,j)=2c_{j-i}$.

\medskip
We also use Brascamp-Lieb inequality in the following formulation.
Let us introduce the following measure,
\begin{equation}
\mu^{0,t}_\Lgr(\dd\ul h) = \frac 1{Z^{0,t}_\Lgr} e^{t\sprod\agr h}
\mu^0_\Lgr(\dd\ul h)\,.
\end{equation}
Expectation value and variance w.r.t. $\mu^{0,t}_\thebox$ are written
$\bk{\,\cdot\,}^{0,t}_\thebox$ and $\var^{0,t}_\thebox$. Then for any
$\agr:\Lgr\ra\bbR$; then
\begin{equation}\label{eq_BL}
\var_\Lgr^{0,t} (\sum_{i\in\Lgr} \sprod \agr h) \leq \frac 1c \var_\Lgr^{0,{\rm G}}
(\sum_{i\in\Lgr} \sprod \agr h)\,,
\end{equation}
where $\var_\Lgr^{0,{\rm G}}$ is the variance w.r.t. the Gaussian measure with
$0$-b.c. in $\Lgr$, obtained by setting $\Cpl_k(x)=x^2/2$ if $P_c(0,k)=1$ and $0$
otherwise (see beginning of the section).

\bigskip\noindent
{\bf Acknowledgments:} The authors thank Volker Bach, Erwin Bolthausen and Dima
Ioffe for interesting discussions on these topics. They also thank Erwin
Bolthausen for communicating the work \cite{BB} before publication.
\section{Mean square height of the pinned field and tail estimates}\label{s_pinV}
This section is devoted to the proof of Theorems \ref{thm_pinV} and
\ref{thm_tail}.
\subsection{Mean square}
We prove now Theorem \ref{thm_pinV}.

Expectation value with respect to $\mu^\pot_\thebox$ has the following
convenient representation, close to the one used in \cite{BI} and
\cite{BB} in the case of the $\dgr$-pinning,
\begin{align}
\bk{\,\cdot\,}^\pot_\thebox
&= \frac 1{Z^\pot_\thebox} \int\dd\ul{h} \,\cdot\,
e^{-\sum_{\rnb{ij}\cap\Lgr\neq\emptyset} \Cpl_{j-i}(h_i-h_j)}
\prod_{j\in\thebox} \Bigl( 1+(e^\egr - 1)\cf{\abs{h_j}\leq a} \Bigr)\nnm\\
&= \sum_{A\subset\thebox} (e^\egr - 1)^{\abs A} \frac
{Z_\thebox(A)}{Z^\pot_\thebox}\; \bk{\,\cdot\;|\,\abs{h_j}\leq a,\,\forall j\in
A}^0_\thebox\,,
\label{eq_repr}
\end{align}
where $Z_\thebox(A) \defby \int\dd\ul{h}\;
e^{-\sum_{\rnb{ij}\cap\Lgr\neq\emptyset} \Cpl_{j-i}(h_i-h_j)}
\prod_{j\in A}\cf{\abs{h_j}\leq a}$.

An upper bound on the mean square height of the field is easily obtained using
\req{eq_repr}. Indeed, we can write
\begin{equation}
\bk{h_i^2}^\pot_\thebox = \sum_{A\subset\thebox} (e^\egr - 1)^{\abs A} \frac
{Z_\thebox(A)}{Z^\pot_\thebox}\; \bk{h_i^2\,|\,\abs{h_j}\leq a,\,\forall j\in
A}^0_\thebox\,.
\end{equation}
Using Lemma \ref{lem_truezero} and \req{eq_BL}, we get
\begin{equation}
\bk{h_i^2\,|\,\abs{h_j}\leq a,\,\forall j\in A}^0_\thebox \leq
4a^2+4\bk{h_i^2}^0_{\Acomp} \leq 4a^2+\frac4c \bk{h_i^2}^{0,{\rm
G}}_{\Acomp}\,,
\end{equation}
where $\Acomp \defby \thebox\setminus A$.

Observe that the random-walk representation of Section \ref{s_statements} gives
\begin{equation}
\bk{h_i^2}^{0,{\rm G}}_{\Acomp} = \ERW{i}{\int_0^\infty\cf{\eta_s=i}\cf{T_A>s}\dd
s}\,,
\end{equation}
where $\ERW{i}{\,\cdot\,}$ denotes expectation with respect to the random walk
starting at the site $i$, $\eta_s$ is the position of the RW at time $s$, and
$T_A\defby\inf_{s\geq 0}\{\eta_s\in A\}$. This last expression can be easily
bounded using a well-known result about symmetric, irreducible random walks (see
e.g. {\bf P}12.3 in \cite{Spitzer}); we obtain
\begin{equation}\label{eq_varest}
\bk{h_i^2}^{0,{\rm G}}_{\Acomp} \leq \frac {\widetilde C}c \log \dd(i,A)\,,
\end{equation}
for some absolute constant $\widetilde C$.
Let $R_{\rm min}$ be the smallest value of the diameter of sets $B$ for which
Proposition \ref{pro_drydecay} applies. Since the range of the random-walk is
$r$-connected, we can use our main estimate, Proposition \ref{pro_drydecay},
which shows that there exists $K>0$ such that ($B_R(i)$ is the ball with radius
$R$ and center $i$)
\begin{align}
\bk{h_i^2}^\pot_\thebox
&\leq 4 a^2 + \frac Cc\log R_{\rm min} + \sum_{R\geq R_{\rm min}}
\sumthree{A\subset\thebox}{A\cap B_R(i)=\emptyset}{A\cap
B_{R+1}(i)\neq\emptyset} \frac Cc \log R\nnm\\
&\leq 4 a^2 + \frac Cc\log R_{\rm min} + \sum_{R\geq R_{\rm min}}
e^{-K\;R^2}\;\frac Cc \log R\,.
\end{align}
This completes the proof of Theorem \ref{thm_pinV}; the estimate on $C_3$
follows by taking the optimal $R_{\rm min}$ above.
\subsection{Tail estimate}
We prove now Theorem \ref{thm_tail}. This proof is close to the previous
one. Let us first prove the upper bound.
Using the representation \req{eq_repr}, we can write
\begin{align}
\bk{\cf{h_i>T}}^V_\thebox
&= \sum_{A\subset\thebox} \nu(A)\; \bk{\cf{h_i>T} \given \abs{h_j}\leq a,\,
\forall j\in A}^0_\thebox \nnm\\
&= \sum_{R\geq 1} \sumthree{A\subset\thebox}{A\cap B_R(i)=\emptyset}{A\cap
B_{R+1}(i)\neq \emptyset} \nu(A) \bk{\cf{h_i>T} \given \abs{h_j}\leq a,\,
\forall j\in A}^0_\thebox\,,
\label{eq_maintail}
\end{align}
where $\nu(A) = (e^\egr-1)^{\abs A} Z_\thebox(A)\big/Z^V_\thebox$. Lemma
\ref{lem_probtruezero} gives
\begin{equation}
\bk{\cf{h_i>T} \given \abs{h_j}\leq a,\,
\forall j\in A}^0_\thebox \leq \bk{\cf{h_i>T-a}}^0_{\Acomp}\,.
\end{equation}
Now this probability is easily evaluated: There exists $\ul C>0$
such that
\begin{equation}\label{eq_dev}
\bk{\cf{h_i>T-a}}^0_{\Acomp} \leq \exp(-\ul C\;T^2/\log R)\,.
\end{equation}
Indeed, this follows from Chebyshev's inequality, Brascamp-Lieb inequality
\req{eq_BL} and the variance estimate \req{eq_varest}.

Let $R_{\rm min}$ be large enough so that we can apply our main estimate to sets
$B$ with $\diam B>R_{\rm min}$. We get
\begin{equation}
\bk{\cf{h_i>T}}^V_\thebox \leq e^{-\cO(T^2)}+\sum_{R\geq R_{\rm min}} e^{-\ul K
R^2 - \ul C \frac {T^2}{\log R}}\,.
\end{equation}
We now have to find the asymptotic behaviour in $T$ of this sum.
Observe that the function $K R^2 + \ul C \frac {T^2}{\log R}$ is convex,
with a unique minimum at $R_0$ solution of
\begin{equation}
R_0\log R_0 = \sqrt{\ul C/2K}\;T\,.
\end{equation}
We cannot solve this equation, however we can easily find a lower bound
on $R_0$:
\begin{equation}
R_0> \sqrt{\ul C/2K}\;\frac T{\log T} \equiv \ol R\,.
\end{equation}
Observe that
\begin{equation}
1\geq \frac {\ol R}{R_0} \geq 1-\cO(\frac{\log\log T}{\log T})\,.
\end{equation}
The required upper bound is obtained by splitting the sum in the following way:
\begin{equation}
\sum_{R\geq 1} e^{-\ul K R^2 - \ul C \frac {T^2}{\log R}}
= \sum_{R=1}^{T} e^{-\ul K R^2 - \ul C \frac {T^2}{\log R}} +
\sum_{R>T} e^{-\ul K R^2 - \ul C \frac {T^2}{\log R}}\,.
\end{equation}
The exponential in the first sum is maximum when $R=R_0$. Therefore,
\begin{equation}
\sum_{R=1}^{T} e^{-\ul K R^2 - \ul C \frac {T^2}{\log R}} \leq 
T\; e^{-\ul K R_0^2 - \ul C \frac {T^2}{\log R_0}} \leq 
T\; e^{-\ul C\;\frac{T^2}{\log T}(1+\cO(\frac{\log\log T}{\log T}))}\,.
\end{equation}
The other part of the sum is easily taken care of by using the bound
\begin{equation}
e^{-\ul K R^2 - \ul C \frac {T^2}{\log R}} \leq e^{-\ul K R^2}\,,
\end{equation}
and estimating the corresponding sum. This finally proves that
\begin{equation}\label{upper}
\sum_{R\geq 1} e^{-\ul K R^2 - \ul C \frac {T^2}{\log R}}
\leq e^{-\ul C\;\frac{T^2}{\log T}(1+\cO(\frac{\log\log T}{\log T})}\,.
\end{equation}

\bigskip
The proof of the lower bound is very similar.
The main change is that we have to use some kind of reverse Chebyshev's
inequality to bound $\bk{\cf{h_i>T+a}}^0_{\Acomp}$. This can be done
by using the following well-known inequality (see \cite{DS} for example),
\begin{equation}
\log\frac {\mu^0_{\Acomp}(h_i>T+a)}{\mu^{0,\agr}_{\Acomp}(h_i>T+a)} \geq -
\frac
{H(\mu^{0,\agr}_{\Acomp}\,|\,\mu^0_{\Acomp})+e^{-1}}{\mu^{0,\agr}_{\Acomp}
(h_i>T+a)}\,,
\end{equation}
where $H(\mu\,|\,\nu)$ is the relative entropy of $\mu$ w.r.t. $\nu$, and
\begin{equation}
\mu^{0,\agr}_{\Acomp}(\dd\ul h) = \frac 1{Z^{0,\agr}_{\Acomp}} e^{\agr h_i}
\mu^0_{\Acomp}(\dd\ul h)\,.
\end{equation}
Differentiating $\bk{h_i}^{0,\agr}_{\Acomp}$ and using \req{eq_BL} and the
reverse Brascamp-Lieb inequality \cite{DGI} to bound the resulting variance in
terms of the corresponding Gaussian quantity, we easily get
\begin{equation}
(Cc)^{-1}\,\agr\,\log R \geq c^{-1}\,\agr\,\var^{\rm G}_{\Acomp}(h_i) \geq
\bk{h_i}^{0,\agr}_{\Acomp} \geq c\, \agr\, \var^{\rm G}_{\Acomp}(h_i) \geq Cc
\,\agr\,\log R\,,
\end{equation}
for some $C>0$; the last inequality follows from well-known result on symmetric,
irreducible random walks, as above. This implies that
\begin{equation}
H(\mu^{0,\agr}_{\Acomp}\,|\,\mu^0_{\Acomp}) \leq \agr
\bk{h_i}^{0,\agr}_{\Acomp} \leq (C c)^{-1} \agr^2 \log R\,. 
\end{equation}
We also have
\begin{equation}
\mu^{0,\agr}_{\Acomp}(h_i>T+a) \geq 1- e^{\agr(T+a)}e^{-\tfrac12 Cc\agr^2\log R} = 1-
e^{-\agr(\tfrac12 \agr C c \log R - T - a)}\,.
\end{equation}
Choosing $\agr=4(T+a)\big/Cc\log R$, this yields
\begin{equation}\label{eq_revcheb}
\mu^0_{\Acomp}(h_i>T+a) \geq (1-e^{-\frac {4(T+a)^2}{C c \log R}})\exp\{-\frac
{4(T+a)^2}{Cc\log R}\Big/(1-e^{-\frac {4(T+a)^2}{C c \log R}})\}\,.
\end{equation}
To obtain the desired lower bound, it suffices to restrict the sum over $R$ in
\req{eq_maintail} to the single term $R=\ol R$, apply Lemma
\ref{lem_probtruezero} to replace $\bk{\cf{h_i>T} \given \abs{h_j}\leq a,\,
\forall j\in A}^0_\thebox$ by $\bk{\cf{h_i>T+a}}^0_{\Acomp}$, and finally use
\req{eq_revcheb}.
\section{Results for the $\dgr$-pinning}\label{s_expdecay}
In this section, we prove Theorems \ref{thm_sgap}, \ref{thm_expdecayGauss} and
\ref{thm_expdecaylargeJ}, and Corollary \ref{cor_odd}.

\medskip
\subsection{Estimate of variances}
We first prove Theorem \ref{thm_sgap} and its corollary. The proof is quite
similar to those given in the previous section. We introduce the symmetric
operator
\begin{equation}
\frK_N (i,j) = \ERW{i}{\int_0^{\tau_{B_N(i)}}\cf{\eta_s=j}\dd s}\,,
\end{equation}
where $B_N(i)=\setof{j\in\bbZZ}{\normI{j-i}<N}$ and
$\tau_{B_N(i)}=\inf\setof{s\geq 0}{\eta_s\not\in B_N(i)}$. From an easy
adaptation of (1.21) in \cite{Lawler}, we know that there exists a constant $C$
such that
\begin{equation}
\sup_i \sum_j \frK_N (i,j) = \sup_i \ERW{i}{\tau_{B_N(i)}} \leq C N^2\,.
\end{equation}
Therefore, using \req{eq_BL} and our main estimate,
\begin{align}
\var_\thebox^J (\sum_i\agr(i)h_i)
&= \sum_{A\subset\thebox} e^{J\abs A} \frac{Z^0_{\Acomp}}{Z^J_\thebox}
\var_{\Acomp}^0 (\sum_i\agr(i)h_i)\nnm\\
&\leq \frac 1c \sum_{A\subset\thebox} e^{J\abs A}
\frac{Z^0_{\Acomp}}{Z^J_\thebox} \var_{\Acomp}^{0,{\rm G}} (\sum_i\agr(i)h_i) \nnm\\
&\leq \frac 1c \sum_{i\in\thebox} \abs{\agr(i)} \ERW{i}{\int_0^{\tau_{B_M(i)}}
\abs{\agr(\eta_s)}\dd s}\nnm\\
&\hspace{1cm} + \frac1c \sum_{i\in\thebox} \abs{\agr(i)} \sum_{N\geq M}
\ERW{i}{\int_{\tau_{B_N(i)}}^{\tau_{B_{N+1}(i)}} \abs{\agr(\eta_s)}\dd s}e^{-KN}
\nnm\\
&\leq \frac 1c \sprod{\abs\agr}{\frK_M\abs\agr} + \frac1c \sum_{N\geq M}
\sprod{\abs\agr}{\frK_{N+1}\abs\agr} e^{-KN}\nnm\\
&\leq \ol C \sum_{i\in\thebox}\agr(i)^2\,.
\end{align}
provided we choose $M$ large enough. This proves Theorem \ref{thm_sgap}.

\bigskip
Let us prove now Corollary \ref{cor_odd}.  We have, similarly as above,
\begin{align}
\var^J_\thebox(F) &= \sum_{A\subset\thebox} \nu(A) \bk{F^2}^0_{\Acomp} - \Bigl\{
\sum_{A\subset\thebox} \nu(A) \bk{F}^0_{\Acomp} 
\Bigr\}^2 \nnm\\
&= \sum_{A\subset\thebox} \nu(A) \var^0_{\Acomp}(F) +
\sum_{A\subset\thebox}  \nu(A) \bigl( \bk{F}^0_{\Acomp} \bigr)^2 -
\Bigl\{ \sum_{A\subset\thebox} \nu(A) \bk{F}^0_{\Acomp} 
\Bigr\}^2\,.
\end{align}
Since $F$ is odd in each coordinate, we know that
\begin{equation}
\bk{F}^0_{\Acomp} = 0\quad\quad\forall A\subset\thebox\,.
\end{equation}
Also, by a version of Brascamp-Lieb inequality proved in \cite{DGI},
\begin{equation}
\var^0_{\Acomp}(F) \leq \frac 1c \var^{0,{\rm G}}_{\Acomp}(\widehat F)\,,
\end{equation}
where $\widehat F(h) \equiv \sum_i \norminf{f_i} h_i$. Now the result follows
from Theorem \ref{thm_sgap}.
\subsection{Mass generation}
The proofs of the last two theorems follow closely the approach of \cite{BB}; they are
based on the representation \req{eq_repr}, which yields the following
expression for the 2-point function (valid for any interactions),
\begin{align}
\bk{h_i h_j}^J_\thebox
&= \sum_{A\subset\thebox} e^{J \abs A} \frac {Z^0_{\Acomp}}{Z^J_\thebox}\;
\bk{h_i h_j}^0_{\Acomp}\nnm\\
&=\sum_{A\subset\thebox} e^{J \abs A} \frac {Z^0_{\Acomp}}{Z^J_\thebox}\;
\ERWNG{i}{A}{\int_0^\infty \cf{\eta_s=j}\cf{T_A>s} \dd
s}\,,
\label{eq_repr2ptf}
\end{align}
where $\ERWNG{i}{A}{\,\cdot\,}$ denotes expectation value w.r.t. the random
walk in random environment described at the end of Section \ref{s_statements}.
Let us first prove Theorem \ref{thm_expdecayGauss}. In this case, the
expectation with respect to the random walk is independent of $A$ and therefore
can be permuted with the sum over $A$, as was done in Section \ref{s_pinV}.
Without loss of generality, we restrict our attention to sites $i$ and $j$
which are sufficiently far from one another so that we can use our main
estimate, Proposition \ref{pro_drydecay}, to get
\begin{equation}
\bk{h_i h_j}^J_\thebox = \int_0^\infty
\ERW{i}{\cf{\eta_s=j}\sumtwo{A\subset\thebox}{A\cap\eta_{[0,s]}=\emptyset}
e^{J\abs A} \frac{Z^0_{\Acomp}}{Z^J_\thebox}}\dd s \leq \int_0^\infty
\ERW{i}{\cf{\eta_s=j}e^{-K\abs{\eta_{[0,s]}}}}\dd s\,.
\end{equation}
This can be estimated as in \cite{BB}. We give here a proof for completeness.
Writing
\begin{equation}
G_N(i,j) \defby \int_0^\infty \ERW{i}{\cf{\eta_s=j}\cf{\tau^{i}_N>s}}\dd s\,,
\end{equation}
with $\tau^{i}_N\defby\inf_{s\geq 0}\{\norminf{\eta_s-i}\geq N\}$, we get
\begin{align}
\bk{h_i h_j}^J_\thebox
&\leq \ERW{i}{\int_0^\infty \cf{\eta_s=j} e^{-K\abs{\eta_{[0,s]}}}\dd s}\nnm\\
&\leq \sum_{N\geq \normI{j-i}}(G_{N+1}(i,j)-G_N(i,j)) e^{-K N}\nnm\\
&\leq (1-e^{-K})\; \sum_{N\geq \normI{j-i}}G_N(i,j)e^{-K(N-1)}\nnm\\
&\leq (1-e^{-K})\; \sum_{N\geq \normI{j-i}}G_N(i,i)e^{-K(N-1)}\nnm\\
&\leq e^{-C_6 \normI{i-j}}\,.
\end{align}
This concludes the proof of Theorem \ref{thm_expdecayGauss}. We don't know how
to make the corresponding  proof in the non-Gaussian case, since in that case
the expectation with respect to the random walk {\em does} depend on the set
$A$, so it is not possible to permute it with the sum over sets $A$. Moreover,
we cannot use the trick of Section \ref{s_pinV} to recover the Gaussian case,
since Brascamp-Lieb inequality does not apply to
$\bk{h_ih_j}^0_{\Acomp}$. It is however possible to do something
when $J$ is large enough, as shown now. Let us write $i \AcConn j$ if $i$ and
$j$ belong to the same $r$-connected (see Section \ref{s_mainestimate})
component of $\Acomp$. Since $\bk{h_i h_j}^0_{\Acomp}
\neq 0$ only if $i \AcConn j$, we can write
\begin{align}
\bk{h_i h_j}^J_\thebox
&= \sumtwo{A\subset\thebox}{i\AcConn j} e^{J \abs A} \frac
{Z^0_{\Acomp}}{Z^J_\thebox}\; \bk{h_i h_j}^0_{\Acomp}\,.
\nnm\\
&\leq \tfrac12 \sumtwo{A\subset\thebox}{i\AcConn j} e^{J \abs A} \frac
{Z^0_{\Acomp}}{Z^J_\thebox}\; (\bk{h_i^2}^0_{\Acomp}+\bk{h_j^2}^0_{\Acomp})\nnm\\
&\leq \tfrac12 \sumtwo{A\subset\thebox}{i\AcConn j} e^{J \abs A} \frac
{Z^0_{\Acomp}}{Z^J_\thebox}\; (\bk{h_i^2}^{0,{\rm
G}}_{\Acomp}+\bk{h_j^2}^{0,{\rm G}}_{\Acomp})\,.
\end{align}
Each of these two terms can be decomposed in the following way.
\begin{align}
\sumtwo{A\subset\thebox}{i\AcConn j} e^{J \abs A} \frac
{Z^0_{\Acomp}}{Z^J_\thebox}\; \bk{h_i^2}^{0,{\rm G}}_{\Acomp}
&\leq \back\int_0^\infty \ERW{i}{\cf{\eta_s=i}\cf{\abs{\eta_{[0,s]}}\geq\normI{i-j}}
\back\sumtwo{A\subset\thebox}{A\cap\eta_{[0,s]}=\emptyset}\back e^{J\abs A}
\frac{Z^0_{\Acomp}}{Z^J_\thebox}\dd s}\nnm\\
&\;+ \int_0^\infty \ERW{i}{\cf{\eta_s=i}\cf{\abs{\eta_{[0,s]}}<\normI{i-j}}
\sumtwo{A\subset\thebox}{i\AcConn j} e^{J\abs A}
\frac{Z^0_{\Acomp}}{Z^J_\thebox}\dd s}\,.
\label{eq_2terms}
\end{align}
The first of these integral can be dealt with as before. To control the second
one, observe that our main estimate Proposition \ref{pro_drydecay} implies the
existence of a constant $\widehat K>0$ such that
and
\begin{equation}
\sumtwo{A\subset\thebox}{i\AcConn j} e^{J\abs A}
\frac{Z^0_{\Acomp}}{Z^J_\thebox} \leq \sumtwo{B\subset\thebox}{i
\BConn j} \sumtwo{A\subset\thebox}{A\cap B=\emptyset} e^{J\abs A}
\frac{Z^0_{\Acomp}}{Z^J_\thebox} \leq \sumtwo{B\subset\thebox}{i
\BConn j} e^{-K\abs B}
\leq e^{-\widehat K
\normI{j-i}}\,,
\end{equation}
as soon as $K$ is large enough (which is true if $J$ is sufficiently large);
here the sum is over sets $B$ which are $r$-connected and contain $i$ and $j$.
Using this, the second integral in \req{eq_2terms} can easily be seen to decay
exponentially with $\normI{i-j}$; Theorem \ref{thm_expdecaylargeJ} follows.
\section{Proof of the main estimate}\label{s_mainestimate}
This section is devoted to the proof of Proposition \ref{pro_drydecay}, which
is the main estimate of this paper. The most important is the first statement,
but the other also appear to be useful. This proposition roughly states that an
arbitrarily weak pinning potential is sufficient to decrease (strictly) the
free energy; its power, however, lies in the fact that it is not restricted to
well-behaved subsets (in the sense of Van Hove for example), but even applies to
``one-dimensional'' ones. 

We say that a set $D\subset\bbZZ$ is {\em $M$-connected}, if, for any $x,y\in
D$, there exists an ordered sequence $(t_0\equiv x,t_1,\dots,t_n\equiv y)$ of
sites of $D$ such that $\normI{t_k-t_{k-1}}\leq M$, for all $k=1,\dots, n$. The
{\em diameter} of a set $D$ is defined by $\diam D=\max_{x,y\in D}\normI{x-y}$.
\begin{pro}\label{pro_drydecay}\hfill

\noindent
1. Let $B\subset\thebox$ be $M$-connected and such that $\diam B \geq
\bigl(a(e^\egr-1)\sqrt{c}\bigr)^{-C(M)}$ for some
$C(M)$ large enough. Then, there exists $K=K(a(e^\egr-1)\sqrt c,M)$, independent of
$B$, such that
$$
\sumtwo{A\subset\thebox}{A\cap B=\emptyset}\; (e^\egr-1)^{\abs A}\;
\frac{Z_\thebox(A)}{Z^V_\thebox} \leq \exp\{-K\,\abs B\}\,.
$$
Moreover, if $a(e^\egr-1)\sqrt c$ is small enough, then there exists $C_8=C_8(M)$ such
that $K>(a(e^\egr-1)\sqrt c)^{C_8}$.\\
2. For any $B\subset\thebox$,
$$
\sumtwo{A\subset\thebox}{A\cap B=\emptyset}\; (e^\egr-1)^{\abs A}\;
\frac{Z_\thebox(A)}{Z^V_\thebox} \geq \exp\{-\egr\,\abs B\}\,.
$$
3. For all $\xi<1$, there exists $C_9=C_9(a(e^\egr-1)\sqrt{c},M)$ such that, for
all $M$-connected $B\subset\thebox$ and all $L$,
$$
\sumtwo{A\subset\thebox}{A\supset B} (e^\egr-1)^{\abs A}
\frac{Z_\thebox(A)}{Z^V_\thebox} \geq (1-\xi)\bigl( 1+ (C_9a(e^\egr-1)\sqrt
c)^{-1} \bigr)^{-\abs B}\,. 
$$
\end{pro}
{\bf Remark:} In the case of the $\dgr$ potential, a look at the proof shows
that the statement corresponding to point 2. takes the form
\begin{equation}
\sumtwo{A\subset\thebox}{A\cap B=\emptyset}\; e^{J\abs A}\;
\frac{Z^0_{\Acomp}}{Z^J_\thebox} \geq (1+e^J)^{-\abs B}\,.
\end{equation}
\begin{proof}
Let us first 1. We introduce the following notations:
\begin{align}
D^k &\defby \{t\in\thebox\,:\,\dd_1(t,D)\leq k\}\,,\\
\bext D \defby D^1\setminus D\,.
\end{align}
The weights $\ol\rho(A) \defby (e^\egr-1)^{\abs
A}\;\frac{Z_\thebox(A)}{Z^V_\thebox}$ define a probability measure on
$\{A\subset\thebox\}$, which we denote by $\Prob$. We also use the notation
$\ol A \defby A\cup\bext\thebox$. What we want to obtain is a upper bound on
\begin{align}
\Prob[A\cap B =\emptyset]
&= \sum_{k\geq 0} \Prob[A\cap B^k =\emptyset\text{ and }\ol A\cap\bext B^k
\neq\emptyset]\nonumber\\
&\leq \sum_{k\geq 0} \Prob[A\cap B^k =\emptyset\,|\,\ol A\cap\bext B^k
\neq\emptyset]\,.
\end{align}
Observe that $B^k$ is also $M$-connected and $\diam B^k>\diam B$. We can write
\begin{align}
\Prob[A\cap B^k =\emptyset\,|\,\ol A\cap\bext B^k \neq\emptyset]
&= \frac {\ds\sumtwo{A\subset \thebox\setminus B^k}{\ol A\cap\bext
B^k\neq\emptyset} (e^\egr-1)^{\abs A}Z_\thebox(A)} {\ds\sumtwo{A\subset
\thebox}{\ol A\cap\bext B^k\neq\emptyset} (e^\egr-1)^{\abs
A}Z_\thebox(A)}\nonumber\\
&= \frac {\ds\sumtwo{A\subset \thebox\setminus
B^k}{\ol A\cap\bext B^k\neq\emptyset} (e^\egr-1)^{\abs{A}}Z_\thebox(A)}
{\ds\sum_{C\subset B^k}(e^\egr-1)^{\abs C}\sumtwo{A\subset \thebox\setminus
B^k}{\ol A\cap\bext B^k\neq\emptyset} (e^\egr-1)^{\abs A}Z_\thebox(A\cup
C)}\nonumber\\
&= \Biggl\{ \sum_{C\subset B^k} (e^\egr-1)^{\abs C} \sumtwo{A\subset \thebox\setminus
B^k}{\ol A\cap\bext B^k\neq\emptyset} \rho(A) \frac{Z_\thebox(A\cup
C)}{Z_\thebox(A)}\Biggr\}^{-1}\nonumber\\
&\leq \Biggl\{ \sum_{C\subset B^k} (e^\egr-1)^{\abs C} \inftwo{A\subset \thebox\setminus
B^k}{\ol A\cap\bext B^k\neq\emptyset} \frac{Z_\thebox(A\cup
C)}{Z_\thebox(A)}\Biggr\}^{-1}\,,
\label{eq_Bdry}
\end{align}
where $\rho(A) \defby (e^\egr-1)^{\abs A} Z_\thebox(A) \big/ \sumtwo{A\subset
\thebox\setminus B^k}{\ol A\cap\bext B^k\neq\emptyset} (e^\egr-1)^{\abs A}
Z_\thebox(A)$.

\begin{figure}[t]
\epsfysize=65mm
\centerline{\epsffile{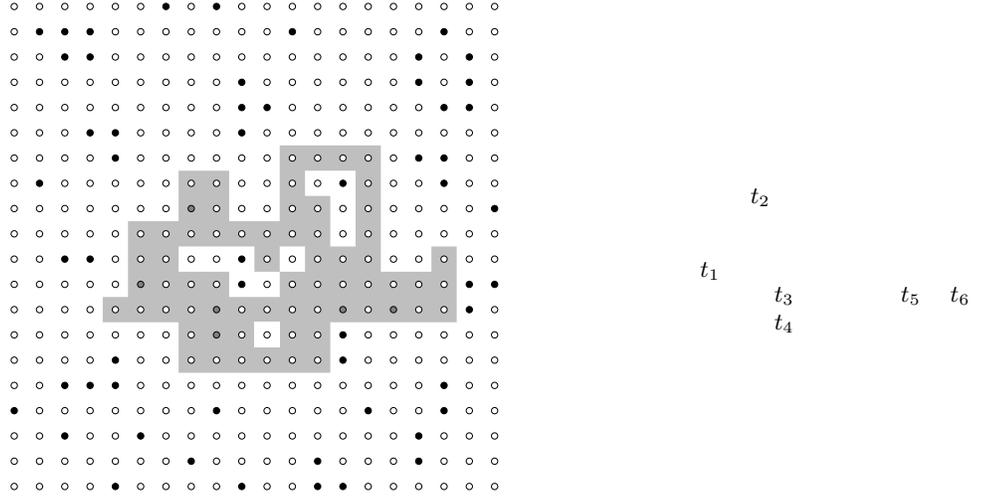}}
\figtext{ 
\writefig	6.92	2.97	{\footnotesize $t_3$} 
\writefig	6.92	2.6	{\footnotesize $t_4$} 
\writefig	6.6	4.28	{\footnotesize $t_2$} 
\writefig	5.93	3.3	{\footnotesize $t_1$} 
\writefig	8.6	2.97	{\footnotesize $t_5$} 
\writefig	9.26	2.97	{\footnotesize $t_6$} 
}
\caption{The set $A$ is represented by the black dots and the set
$C=\{t_1,\dots,t_6\}$ by the
grey ones; the set $B^k$ (connected, here) is composed
of all the sites centered on a shaded plaquette. The set $A_k$ is composed of
the union of $A$ and the sites $t_1,\dots,t_k$.} 
\label{fig_ABC}
\end{figure} 
One has therefore to bound the ratio of partition functions. If we enumerate
the elements of $C$, say $C=\{t_1,\dots,t_{\abs C}\}$, and define
$A_k \defby A \cup \{t_1,\dots,t_k\}$, we get
\begin{equation}
\frac{Z_\thebox(A\cup C)}{Z_\thebox(A)} =
\frac{Z_\thebox(A_1)}{Z_\thebox(A_0)}
\frac{Z_\thebox(A_2)}{Z_\thebox(A_1)}\cdots
\frac{Z_\thebox(A_{\abs C})}{Z_\thebox(A_{\abs C -1})}\,.
\end{equation}
But, using Lemmas \ref{lem_probtozero} and \ref{lem_lowerbd},
\begin{equation}\label{eq_ratio}
\frac{Z_\thebox(A_k)}{Z_\thebox(A_{k-1})} =
\mu^0_\thebox(\abs{h_{t_k}}\leq a \given \abs{h_j}\leq a,\,\forall j\in A_{k-1}) \geq
\tfrac12 \mu^0_{\Acomp_{k-1}}(\abs{h_{t_k}}\leq a) \geq \frac a
{8\bk{\abs{h_{t_k}}}^0_{\Acomp_{k-1}}}\,.
\end{equation}
Therefore,
\begin{equation}
\frac{Z_\thebox(A\cup C)}{Z_\thebox(A)} \geq \prod_{k=1}^{\abs
C} \frac a {8\bk{\abs{h_{t_k}}}^0_{\Acomp_{t_{k-1}}}}\,.
\end{equation}
To go further, we need to use the properties that $B^k$ inherited from $B$. 
Let $l\in {\mathbb N}$ large enough (in particular $l \gg M$), but small
compared to $\diam B$; we
consider a grid of spacing $l$ in $\Lambda$, with cells $\cC_i$.
Observe that there exists two numbers $\nu\in(0,1]$ and $\rho\in(0,1]$,
independent of the set $B^k$ and of $l$, such that the following properties hold:
\begin{itemize}
\item $B^k\subset \bigcup_{j\in\bbJ} \cC_j$,
\item $B^k\cap\cC_j\neq\emptyset$, for all $j\in\bbJ$,
\item $\abs{B^k\cap\cC_j}>\frac\nu M l$, for all $j\in\tilde\bbJ$,
\end{itemize}
where $\{\cC_j,\,j\in\bbJ\}$ is a connected set of cells, and $\tilde\bbJ
\subset\bbJ$ with $\abs{\tilde\bbJ} \geq \rho \abs{\bbJ}$.
\begin{figure}[t]
\epsfysize=65mm
\centerline{\epsffile{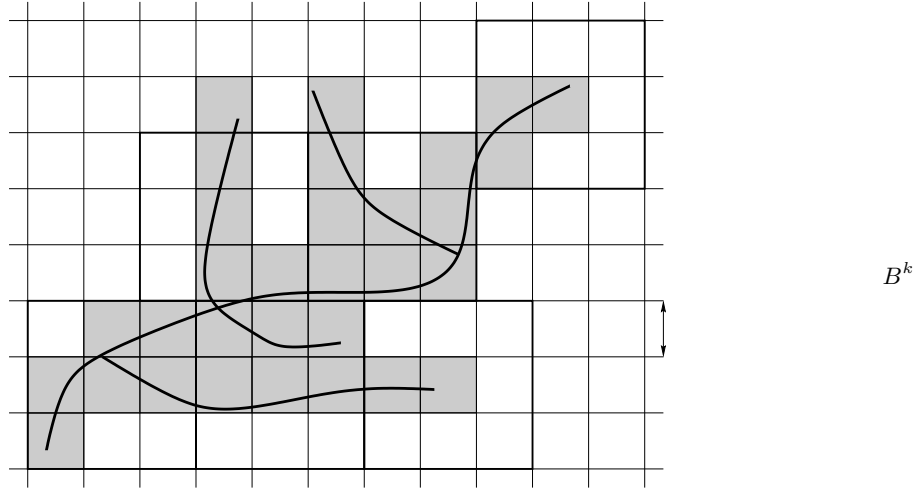}}
\figtext{ 
\writefig	7.27	3.2	{\footnotesize $B^k$} 
\writefig	11.8	2.53	{\footnotesize $l$} 
}
\caption{The shaded cells represents the set $\{\cC_j,\, j\in\bbJ\}$; the six
large $3\times 3$ squares represents the cells $\{\cD_i,\,i=1,\dots,N_D\}$. The
summation will be done on all sets $C\subset B^k$ containing exactly one site in
each of the cells $\cC_j$, $j\in\bbJ$.} 
\label{fig_boxes}
\end{figure} 
Indeed, the first statements are a simple consequence of the $M$-connectedness
of $B^k$, and the last one is proven in the following way. Let
$\{\cD_i,\,i=1,\dots,N_D\}$ be a set of disjoint square boxes in ${\mathbb
Z}^2$, build with exactly 9 cells of the grid defined above, and such that the
middle-cell of each such box belongs to $\{\cC_j,\, j\in\bbJ\}$. We suppose
that these boxes are chosen in such a way as to maximize $N_D$ under these
constraints. Then
\begin{itemize}
\item At most $12 N_D$ cells of $\{\cC_j,\, j\in\bbJ\}$ are outside every
$\cD_i$.
\item Each $\cD_i$ contains at most $9$ cells of $\{\cC_j,\, j\in\bbJ\}$.
\end{itemize}
Therefore, $9 N_D+ 12 N_D \geq \abs \bbJ$, i.e.
\begin{equation}
N_D\geq \tfrac1{21} \abs\bbJ\,.
\end{equation}
Now, $\abs{B^k\cap\cD_i}>l/r$, for all $i=1,\dots,N_D$. Consequently, each box
$\cD_i$ contains at least one cell $\cC$ with $\abs{B^k\cap\cC}\geq \frac 1{9M} l$.
Choosing $\nu=\tfrac19$, this implies that
\begin{equation}
\abs{\tilde{\bbJ}}\geq \tfrac19 N_D \geq \tfrac 1{189} \abs{\bbJ}.
\end{equation}
We can therefore take $\rho=\tfrac1{189}$.

\bigskip
We'll restrict the summation in (\ref{eq_Bdry}) on sets $C\subset B^k$ which
satisfy
\begin{equation}
\abs{C\cap\cC_i}=1\,,\quad\forall i\in \bbJ\,.
\end{equation}
We number the elements of $C$ as above, but in such a way as to ensure that
\begin{equation}
\{\cC\,:\,\cC\ni t_i\,,1\leq i \leq k\}
\end{equation}
is connected for all $1\leq k \leq \abs C = \abs \bbJ$. Then
$\dd_1(t_k,A_{k-1}) \leq \sqrt{5} l$, for all $k>1$. We further ask
that $\dd_1(t_1,\ol A) \leq \sqrt{5} l$, which is always possible. Using
this, we obtain
\begin{equation}
\frac{Z_\thebox(A\cup C)}{Z_\thebox(A)} \geq \left( \frac
{a\ol K\sqrt c}{\sqrt{\log l}} \right)^{\abs\bbJ}\,.
\end{equation}
Indeed, \req{eq_BL} implies
\begin{align}
\bk{\abs{h_{t_k}}}^0_{\Acomp_{k-1}}
&\leq \bigl[ \bk{h_{t_k}^2}^0_{\Acomp_{k-1}} \bigr]^{\tfrac12}\nnm\\
&\leq \bigl[\frac 1c \bk{h_{t_k}^2}^{0,{\rm G}}_{\Acomp_{k-1}}
\bigr]^{\tfrac12}\nnm\\
&\leq \frac {\ol K}{\sqrt c} \sqrt{\log l}\,,
\end{align}
since, by construction, $\dd_1(t_k,A_{t_{k-1}})\leq \sqrt 5 l$, and the
expectation value can be estimated using the random walk representation and
standard results about irreducible, symmetric random walk, see {\bf P}11.6 and
{\bf P}12.3 in \cite{Spitzer} for example. Therefore,
\begin{align}
\sum_{C\subset B^k} (e^\egr-1)^{\abs C} \inftwo{A\subset \thebox\setminus
B^k}{\ol A\cap\bext B^k\neq\emptyset} \frac{Z_\thebox(A\cup C)}{Z_\thebox(A)}
&\geq \sumtwo{C\subset B^k}{\abs{C\cap\cC_i}=1,\,\forall i\in\bbJ} (e^\egr-1)^{\abs\bbJ}
\left( \frac {a\ol K\sqrt c}{\sqrt{\log l}} \right)^{\abs\bbJ}\nonumber\\ 
&\geq (\nu l)^{\rho\abs\bbJ}\left( \frac {(e^\egr-1)a\sqrt c \ol K}{\sqrt{\log l}}
\right)^{\abs\bbJ}\nonumber\\ 
&= \left( \frac {(e^\egr-1)a\sqrt c \ol K(\nu l)^\rho}{\sqrt{\log l}}
\right)^{\abs\bbJ}\nonumber\\
&\geq e^{\abs\bbJ} \geq e^{\abs B /l^2}\,,
\end{align}
which implies
\begin{equation}
\Prob[A\cap B^k =\emptyset\given \ol A\cap\bext B^k \neq\emptyset] \leq
\exp\{-\tilde{K}_l\;\abs {B^k}\}\,,
\end{equation}
for some $\tilde{K}_l>0$ independent of $B^k$ (provided $\diam B \gg l\geq
l_0(a(e^\egr-1)\sqrt c,M)$). Therefore, we finally have
\begin{equation}
\Prob[A\cap B =\emptyset] \leq \exp\{-K\;\abs B\}\,,
\end{equation}
for some $K>0$. The explicit bound on $K$ follows by optimizing over $l_0$
above; this also explains the constraint on  $\diam B$.

\bigskip
Let us prove 2. Proceeding as in \req{eq_Bdry}, we can write
\begin{align}
\sumtwo{A\subset\thebox}{A\cap B=\emptyset}\; (e^\egr-1)^{\abs A}\;
\frac{Z_\thebox(A)}{Z^V_\thebox}
&\geq \Biggl\{ \sum_{C\subset B} (e^\egr-1)^{\abs C} \sup_{A\subset
\thebox\setminus B} \frac{Z_\thebox(A\cup
C)}{Z_\thebox(A)}\Biggr\}^{-1}\nnm\\
&\geq \Biggl\{ \sum_{C\subset B} (e^\egr-1)^{\abs C}\Biggr\}^{-1}\nnm\\
&=e^{-\egr\abs B}\,.
\end{align}

\bigskip
We finally prove 3.
Let $R$ be a large enough integer; we write
\begin{align}
\Prob[A\supset B] &\geq \Prob[A\supset B,\,\ol A\cap(B^R\setminus
B)\neq\emptyset]\nnm\\
&= \Prob[A\supset B \given \ol A\cap(B^R\setminus B)\neq\emptyset]
(1-\Prob[\ol A\cap(B^R\setminus B)=\emptyset])\nnm\\
&\geq (1-\xi) \Prob[A\supset B \given \ol A\cap(B^R\setminus B)\neq\emptyset]\,,
\end{align}
provided $R$ is large enough, by part 1. of the proposition. Now, similarly
as before,
\begin{multline}
\Prob[A\supset B \given \ol A\cap(B^R\setminus B)\neq\emptyset]\\
= \frac{\ds\sumtwo{A\supset B}{\ol A\cap(B^R\setminus B)\neq\emptyset}
(e^\egr-1)^{\abs A}Z_\thebox(A)} {\ds\sumtwo{A\supset B}{\ol A\cap(B^R\setminus
B)\neq\emptyset} (e^\egr-1)^{\abs A}Z_\thebox(A) \Bigl(\sum_{C\subset B}
(e^\egr-1)^{-\abs C} \frac{Z_\thebox(A\setminus C)}{Z_\thebox(A)} \Bigr) }\,.
\end{multline}
The conclusion follows from \req{eq_ratio} and the construction, which imply
that
\begin{equation}
\frac{Z^0(A)}{Z^0(A\setminus C)} \geq  \Bigl( \frac{\sqrt c a\ol K}{\sqrt{\log R}}
\Bigr)^{\abs C}\,.
\end{equation}
Indeed this gives
\begin{align}
\Prob[A\supset B \given \ol A\cap(B^R\setminus B)\neq\emptyset]
&\geq\Bigl\{ \sum_{C\subset B}\bigl( \frac{\sqrt{\log R}}{a(e^\egr-1)\sqrt c\ol
K} \bigr)^{\abs C} \Bigr\}^{-1}\nnm\\
&= \Bigl\{ 1+ \frac{\sqrt{\log R}}{a(e^\egr-1)\sqrt c\ol K}\Bigr\}^{-\abs B}\,.
\end{align}
\end{proof}
\section{Appendix: Proofs of some technical estimates}
In this section, we give the proofs of several technical statements used in the
previous ones. Since FKG inequality is used several times, we recall that, as a
consequence of Corollary 1.7 in \cite{HP}, measures of the form
\begin{equation}
\frac{\mu^b_\Lgr(\,\cdot\,\prod_{i\in\Lgr}f_i(h_i))}
{\mu^b_\Lgr(\prod_{i\in\Lgr}f_i(h_i))}
\end{equation}
are FKG.
\begin{lem}\label{lem_tech}
Let $g$ be a positive, even function which is increasing on $\bbR^+$ and such
that $g(0)=0$. Then, for any $\Lgr\Subset\bbZZ$, any $A\subset\Lgr$ and any
$j\in\Lgr\setminus A$,
$$
\bk{g(h_j)\given \abs{h_k}\leq a\;\forall k\in A}^0_\Lgr \leq
\bk{g(h_j+a)\given h_j\geq -a}^0_{\Lgr\setminus A}\,.
$$
\end{lem}
\begin{proof}
We introduce $g_\nea(h_j) = g(h_j \vee 0)$. Using symmetry, FKG twice and translation
invariance, we can write
\begin{align}
\bk{g(h_j)\given \abs{h_k}\leq a\;\forall k\in A}^0_\Lgr
&= \bk{g_\nea(h_j)\given \abs{h_k}\leq a\;\forall k\in A,\,h_j\geq
0}^0_\Lgr\nnm\\
&\leq  \bk{g_\nea(h_j)\given h_k= a\;\forall k\in A,\,h_j\geq
0}^0_\Lgr\nnm\\
&\leq \bk{g_\nea(h_j)\given h_k= a\;\forall k\in A,\,h_j\geq
0}^a_\Lgr\nnm\\
&= \bk{g_\nea(h_j+a)\given h_k= 0\;\forall k\in A,\,h_j\geq
-a}^0_\Lgr\nnm\\
&= \bk{g(h_j+a)\given h_j\geq -a}^0_{\Lgr\setminus A}\,.
\end{align}
Let us explain how the two inequalities are obtained.
Let $\lgr>0$. Since $\prod_{k\in A}\cf{h_k> a-\lgr}$ and
$g_\nea$ are increasing, and the measure $$\cf{h_j\geq
0} \prod_{k\in A}\cf{\abs{h_k}\leq a} \dd\mu^0_\Lgr$$ is FKG, 
\begin{equation}
\bk{g_\nea(h_j)\,|\,\abs{h_k}\leq a,\,\forall k\in A, h_j\geq 0}^0_\Lgr \leq
\bk{g_\nea(h_j)\,|\,\abs{h_k}\in (a-\lgr, a],\,\forall k\in A, h_j\geq
a}^0_\Lgr\,.
\end{equation}
Letting $\lgr$ go to zero gives the first inequality. The second follows from
the observation that  $\mu^0_\Lgr(\dd \ul h) = \phi_b(\ul h)\mu^b_\Lgr(\dd \ul
h)$, with $\ds\phi_b(\ul h) = \prodtwo{\rnb{ik}}{i\in\Lgr,\,k\not\in\Lgr}
\exp\{\Cpl_{k-i}(h_i-b)-\Cpl_{k-i}(h_i)\}$ decreasing if $b>0$. Indeed,
\begin{equation}
\frac{\dd}{\dd h_i} (\Cpl_{k-i}(h_i-b)-\Cpl_{k-i}(h_i)) = -\int_{h_i-b}^{h_i}
\Cpl_{k-i}''(x)\dd x <0\,, \forall i\in\Lgr\,.
\end{equation}
\end{proof}
\begin{lem}\label{lem_probtruezero}
Let $T>0$. Then, for all $\Lgr\Subset\bbZZ$, $A\subset\Lgr$ and
$i\in\Lgr\setminus A$,
$$
\mu^0_{\Acomp}(h_i>T+a)\leq \mu^0_\Lgr(h_i>T \given \abs{h_j}\leq a,\,\forall
j\in A) \leq \mu^0_{\Acomp}(h_i>T-a)\,.
$$
\end{lem}
\begin{proof}
The proof is completely similar to the previous one. We have
\begin{align}
\bk{\cf{h_i>T} \given \abs{h_j}\leq a\;\forall j\in A}^0_\thebox
&\leq  \bk{\cf{h_i>T} \given \abs{h_j}= a\;\forall j\in A}^0_\thebox \nnm\\
&\leq  \bk{\cf{h_i>T} \given \abs{h_j}= a\;\forall j\in A}^a_\thebox \nnm\\
&\leq  \bk{\cf{h_i>T-a}}^0_{\Acomp} \,,
\end{align}
and
\begin{align}
\bk{\cf{h_i>T} \given \abs{h_j}\leq a\;\forall j\in A}^0_\thebox
&\geq  \bk{\cf{h_i>T} \given \abs{h_j}= -a\;\forall j\in A}^0_\thebox \nnm\\
&\geq  \bk{\cf{h_i>T} \given \abs{h_j}= -a\;\forall j\in A}^{-a}_\thebox \nnm\\
&\geq  \bk{\cf{h_i>T+a}}^0_{\Acomp} \,.
\end{align}
\end{proof}
\begin{lem}\label{lem_truezero}
For any $\Lgr\Subset\bbZZ$ and $A\subset\Lgr$,
$$
\bk{h_i^2\,|\,\abs{h_j}\leq a,\,\forall j\in
A}^0_\Lgr \leq 4a^2+4\bk{h_i^2}^0_{\Lgr\setminus A}\,.
$$
\end{lem}
\begin{proof}
This follows easily from Lemma \ref{lem_tech} with $g(x)=x^2$, and FKG
inequality which yields
\begin{equation}
\bk{((h_i+a) \vee 0)^2\given h_i\geq -a}^0_{\Lgr\setminus A} \leq \bk{((h_i+a)
\vee 0)^2\given h_i\geq 0}^0_{\Lgr\setminus A}\,.
\end{equation}
\end{proof}
\begin{lem}\label{lem_probtozero}
For any $\Lgr\Subset\bbZZ$ and $A\subset\Lgr$,
$$
\mu^0_\Lgr(\abs{h_i}\leq a\,|\,\abs{h_j}\leq a,\,\forall j\in
A)^0_\Lgr \geq \tfrac12 \mu^0_{\Lgr\setminus A}(\abs{h_i}\leq a)\,.
$$
\end{lem}
\begin{proof}
Lemma \ref{lem_tech} with $g(x)=\cf{\abs{x}>a}$ implies that
\begin{align}
\mu^0_\Lgr(\abs{h_i}\leq a\,|\,\abs{h_j}\leq a,\,\forall j\in A)^0_\Lgr
&\geq \mu^0_{\Lgr\setminus A}(\abs{h_i+a}\leq a \given h_i\geq -a)\nnm\\
&= \mu^0_{\Lgr\setminus A}(h_i\leq 0 \given h_i\geq -a)\nnm\\
&\geq \mu^0_{\Lgr\setminus A}(-a\leq h_i\leq 0)\nnm\\
&= \tfrac12 \mu^0_{\Lgr\setminus A}(\abs{h_i}\leq a)\,.
\end{align}
\end{proof}
\begin{lem}\label{lem_lowerbd}
For any $\Lgr\Subset\bbZZ$ and $i\in\Lgr$,
$$
\mu^0_\Lgr(\abs{h_i}\leq a) \geq \frac a{4\bk{\abs{h_i}}^0_\Lgr} \wedge
\tfrac12\,.
$$
\end{lem}
\begin{proof}
The proof follows from the following elementary result, which is proved in
\cite{DMRR}: {\em Let X be a random variable whose density under $\bbP$ is even and
decreasing on $\bbR^+$. Then $\bbP[\abs X \leq a] \geq \frac a{4\bbE[\abs X]}
\wedge \tfrac12$, where $\bbE[\,\cdot\,]$ is the expectation value with respect
to $\bbP$}.

Let $F_j$ be the density of $h_j$ under $\mu^0_\Lgr$. The
evenness is obvious; let us check the monotonicity.
\begin{align}
&F_j(x) = \frac 1{Z^0_\Lgr} \int \back\prod_{t\in\bbR^{\Lgr\setminus\{j\}}}
\back\dd h_t \prodtwo{\rnb{kl}}{k,l\in\Lgr\setminus \{j\}} \back
\e{-\Cpl_{l-k}(h_k-h_l)} \back\prodtwo{\rnb{jk}}{k\in\Lgr}
\back\e{-\Cpl_{k-j}(h_k-x)} \back\prodtwo{\rnb{jk}}{k\not\in\Lgr}
\back\e{-\Cpl_{k-j}(x)} \back\prodthree{\rnb{kl}}{k\in\Lgr\setminus
\{j\}}{l\not\in\Lgr} \back\e{-\Cpl_{l-k}(h_k)}\,,\nnm\\
&F'_j(x) = \tilde C_1 \sumtwo{\rnb{jk}}{k\not\in\Lgr} \bk{-\Cpl_{k-j}'(x)\given
h_j=x}^0_\Lgr + \tilde C_2 \sumtwo{\rnb{jk}}{k\in\Lgr}
\bk{\Cpl_{k-j}'(h_k-x)\given h_j=x}^0_\Lgr\,,
\end{align}
where $\tilde C_1$ and $\tilde C_2$ are two positive constants. Now, for $x\geq
0$, $\Cpl_{k-j}'(x)\geq 0$ and therefore the first term is negative. By FKG,
\begin{equation}
\bk{\Cpl_{k-j}'(h_k-x)\given h_j=x}^0_\Lgr = \bk{\Cpl_{k-j}'(h_k)\given
h_j=0}^{-x}_\Lgr \leq \bk{\Cpl_{k-j}'(h_k)}^0_{\Lgr\setminus
\{j\}}=0\,,
\end{equation}
since $\Cpl_{k-j}'$ is increasing and odd. Consequently, $F'_j(x)\leq 0$ for $x\geq 0$.
\end{proof}

\end{document}